\documentclass[11pt]{amsart}
\usepackage{amssymb,amsfonts,amsmath, amscd, mathrsfs, latexsym, amsthm}
\usepackage[dvips]{graphicx}
\usepackage[all]{xy}
\title{Semi-topological K-theory for certain projective varieties}                 
\theoremstyle{plain} % non-italicized text in theorem environment
\newtheorem{theorem}{Theorem}[section]
\newtheorem{lemma}[theorem]{Lemma} % use the theorem counter
\newtheorem*{definition}{Definition}
\newtheorem*{claim}{Claim}
\newtheorem{conjecture}{Conjecture}[section]
\newtheorem{proposition}{Proposition}[section]
\newtheorem{corollary}{Corollary}[section]
\newtheorem{remark}{Remark}[section]
\newtheorem{convention}{Convention}[section]
\newcommand{\R}{\mathbb{R}}
\newcommand{\C}{\mathbb{C}}
\newcommand{\tensor}{\otimes}
\author{Mircea Voineagu}

\address {Department of Mathematics, Northwestern University,
Evanston, IL  60208}
\email{mircea@math.northwestern.edu}



\begin{document}
%%\author{Mircea Voineagu}

%%\address {Department of Mathematics, Northwestern University,
%%Evanston, IL  60208}
%%% Email address is optional.
%%\email{mircea@math.northwestern.edu}

%%\subjclass[2000]{19E20, 19E15, 14F43}

%%\keywords{}
\begin{abstract}
In this paper we compute Lawson homology groups and semi-topological K-theory for certain threefolds and fourfolds. We consider smooth complex projective varieties whose zero cycles are supported on a proper subvariety. Rationally connected varieties are examples of such varieties. The computation makes use of a technique of Bloch and Srinavas, of the Bloch-Kato conjecture and of the spectral sequence relating morphic cohomology and semi-topological K-theory.

% Put abstract here
\end{abstract}
\maketitle
\tableofcontents
\section{Introduction}
 Eric Friedlander and Mark Walker introduced in \cite{FW} the (singular) semi-topological K-theory of a complex projective variety $X$. This is defined by
$$K^{sst} _*(X)=\pi _*(Mor(X,Grass)^+)$$
where $Grass=\amalg _{n,N}Grass _{n}(\mathbb{P}^N)$. By $Mor(X,Grass)^+$ we denote the topological group given by the homotopy completion of the space of algebraic maps between $X$ and $Grass$.

Semi-topological K-theory lies between algebraic and topological K-theory in the sense that the natural map from the algebraic K-theory $K _{*}(X)$ of a variety $X$ to the connective (complex) topological K-theory $ku^*(X^{an})$ of its underlying analytic space $X^{an}$ factors through the semi-topological K-theory of $X$, i.e.
$$K _q(X)\rightarrow K^{sst} _q(X)\rightarrow ku^{-q}(X^{an})$$
for any $q\geq 0$.

In \cite{L} Blaine Lawson introduced the (Lawson) homology groups of a projective complex variety $X$, which are given by
$$L _rH _n(X)=\pi _{n-2r}(\mathscr{Z} _r(X))$$
where $\mathscr{Z} _r(X)$ is the naive group completion of the topological monoid $$\mathscr{C} _r(X)=\amalg _d\mathscr{C} _{r,d}(X)$$ 
with $\mathscr{C} _{r,d}(X)$ the Chow variety of subvarieties of $X$ of dimension $r$ and degree $d$.

In \cite{FL} Eric Friedlander and Blaine Lawson introduced the morphic cohomology, a cohomology theory dual to Lawson homology \cite{FL1}. They defined 
$$L^rH^n(X)=\pi _{2r-n}(\mathscr{Z}^r(X))$$
where $\mathscr{Z}^r(X)$ is the naive group completion of the following topological monoid $Mor(X,\mathscr{C} _0(\mathbb{P}^r))/Mor(X,\mathscr{C} _0(\mathbb{P}^{r-1})).$

Morphic cohomology groups are related to the semi-topological K-theory by means of a semi-topological spectral sequence \cite{FHW} compatible with the motivic spectral sequence and Atiyah-Hirzebruch spectral sequence \cite{FHW}.

In this paper we study the map 
$$K _*^{sst}(X)\rightarrow ku^{-*}(X^{an})$$
for various complex projective varieties $X$.

We divide the paper in six sections. In the second section we fix the notations and recall some essential results that we need in the paper.

In the third section we study the effects of the Bloch-Kato conjecture on the kernel and cokernel of the generalized cycle maps. We give a new proof of a theorem of Bloch about the torsion of the singular cohomology of a smooth projective variety. We also study the torsion of the Borel-Moore homology of a quasi-projective smooth variety. At the end of the section we construct a birational invariant using Lawson homology.

In the fourth section we study the action of an algebraic cycle on morphic cohomology groups. Our approach is slightly different than the approach pioneered by C. Peters \cite{P} and our results include the results of \cite{P}.

In the fifth section we start comparing Lawson homology and singular homology of smooth projective varieties with zero cycles supported on a subvariety. We essentially use the results of the previous two sections and a technique introduced by Bloch and Srinavas \cite{BS}. 

The main goal of this section is to study the semi-topological K-theory of our ``degenerate'' varieties. One of the main results of the section is the following theorem which computes semi-topological K-theory of ``degenerate''  threefolds.

\begin{theorem}
Let $X$ be a smooth projective complex threefold such that there is a proper subvariety $V\subset X$ with $CH _0(X \setminus V)=0$. Then:
$$K^{sst} _i(X)\simeq ku^{-i}(X^{an}), i\geq 1,$$
$$K^{sst} _0(X)\hookrightarrow ku^0(X^{an}).$$
Moreover if $X$ is a rationally connected threefold then
$$K^{sst} _i(X)\simeq ku^{-i}(X^{an}),i\geq 0.$$
\end{theorem}

This computation generalizes a result of \cite{FHW} about the semi-topological K-theory of a rational threefold. 

The following result describes the semi-topological K-theory of some ``degenerate'' fourfolds.

\begin{theorem}
Let $X$ be a smooth projective fourfold such that there is a proper subvariety $V\subset X$ of $dim(V)\leq 2$ with $CH _0(X \setminus V)=0$. Then: 
$$K^{sst} _i(X)\simeq ku^{-i}(X^{an}),i\geq 3,$$
$$K^{sst} _2(X)\hookrightarrow ku^{-2}(X^{an}),$$
$$K^{sst} _i(X) _\mathbb{Q}\simeq ku^{-i}(X^{an}) _\mathbb{Q},i=1,2,$$
$$K^{sst} _0(X) _\mathbb{Q}\hookrightarrow ku^0(X^{an}) _\mathbb{Q}.$$
\end{theorem}

We may contrast the above results with a result of H. Gillet \cite{G} (see also C. Pedrini and C. Weibel \cite{PW}). He proved that the image of the map 
$$K _n(X)\rightarrow ku^{-n}(X)$$
is finite for any $n>0$ and for any $X$ smooth complex projective variety.

In the sixth section of the paper we give some consequences of a theorem of Jannsen \cite{Ja} and Laterveer \cite{Lat} concerning a special decomposition of the diagonal for varieties with small Chow groups. 

In the last section of the paper we study morphic cohomology of projective smooth linear varieties. The main idea is to use a $K\ddot{u}nneth$ formula for such varieties proved by R. Joshua \cite{RJ} and by B. Totaro \cite{To}. The results in this section were proved in \cite{FHW} using other tools. 

This work started from a question of my advisor Eric Friedlander refering to the semi-topological K-theory of Fano varieties. I gratefully acknowledge his guidance and many valuable suggestions, in particular the suggestion that the Bloch-Kato conjecture may be helpful in the morphic cohomology context. I am also thankful to Mark Walker for carefully reading the manuscript and making useful remarks and to Jeremiah Heller for helpful discussions.
\section{Notations and Recollection}
 Throughout this paper $X$ will denote a smooth projective irreducible variety over the complex numbers of dimension $d$ (unless otherwise stated). By $H^{p,q} _\mathbb{M}(X)$, $L^qH^p(X)$ and $L _pH _q(X)$ we denote motivic cohomology, morphic cohomology and Lawson homology. For a field $E$ we denote $K^M _*(E)$ to be the Milnor K-theory of $E$. By $cyc^{p,q}$, respectively $cyc _{p,q}$ we denote the generalized cycle maps
$$cyc^{p,q}:L^pH^q(X)\rightarrow H^q(X)$$
respectively
$$cyc _{p,q}:L _pH _q(X)\rightarrow H _q(X).$$
We denote by $K^{q,n}:=Ker\{cyc^{q,n}\}$, $K _{p,q}:=Ker(cyc _{p,q})$ and by $C^{q,n}:=Coker\{cyc^{q,n}\},$ $C _{q,n}=Coker\{cyc _{q,n}\}$. For an abelian group $A$ we denote $ _mA:=\{a\in A| ma=0\}$. 

If for a variety $X$ there is a proper subvariety $V\subset X$ such that $CH _0(X\setminus V)=0$ we say as in \cite{BS} that $X$ is ``degenerate''  and also that its zero cycles ``are supported on subvariety $V$''. 

For a complex variety $X$ we denote $X^*$ a resolution of singularities for $X$.

 We will start recalling the basics about the (co)niveau filtration of the singular (co)homology. Let 
$$N _kH _n(X)=\sum _{dim(W)\leq k}Im(H _n(W)\rightarrow H _n(X))$$
be a step in the niveau filtration of $H _n(X)$. This is an ascending filtration$$0\subset N _0H _n(X)\subset ..\subset N _kH _n (X)\subset ...\subset H _n(X)$$
which has the property that 
\begin{equation}
\label{*}
N _kH _n(X)=H _n(X)
\end{equation}
 for any $k\geq min\{n,d\}$.

It is easy to see that $N _dH _n(X)=H _n(X)$ for any natural $n$. For $n<d$ the above equality follows from an induction argument using weak Lefshetz theorem. For $X$ smooth we know that the niveau filtration is isomorphic to the coniveau filtration of the cohomology of $X$, i.e
\begin{equation}
\label{**}
  N _kH _n(X)\simeq N^{d-k}H^{2d-n}(X)
\end{equation}
where we define 
$$N^kH^n(X)=\sum _{cd(W)\geq k}Im(H^n _W(X)\rightarrow H^n(X)).$$
 From (\ref{*}) and from (\ref{**}) we conclude that 
$$N _{d-1}H _{2d-n}(X)\simeq N^1H^n(X)\simeq H^n(X)\simeq H _{2d-n}(X)$$
for any $n$ such that $2d-n\leq d-1\Leftrightarrow n\geq d+1.$
We also know (\cite{FM} and \cite{W}) the following property of the generalized cycle maps 
\begin{proposition}(\cite{FM} and \cite{W})
For a smooth projective variety $X$
$$Im(cyc^{q,n})\subset N^{n-q}H^n(X)$$ 
with equality when $n=2q$ or $n=2q-1$.
\end{proposition}

For a quasi-projective variety, Deligne \cite{D} and Gillet-Soule \cite{GS} defined a weight filtration on the Borel-Moore homology of $U^{an}$ (denoted by $H^{BM} _*(U^{an})$).  We recall the definition of this filtration. Choose a compactification $U\subseteq X$ so that $X$ is a projective complex variety and let $Y$ be the reduced complement of $U$ in $X$. Consider $\mathbb{Z}Sing _*()$ the functor taking a space $Z$ to the complex associated to the simplicial set $Sing _*Z$. We may construct two hypercovers (\cite{GS}) $X _*\rightarrow X$ and $Y _*\rightarrow Y$ such that $X _n$ and $Y _n$ are smooth projective varieties and such that there is a map $Y _*\rightarrow X _*$ which covers the embedding $Y\subset X$. Denoting $U _n=X _n\sqcup Y _{n-1}$ we may construct a bicomplex
\begin{equation}
\label{***}
  ...\rightarrow \mathbb{Z}Sing _*(U _1)\rightarrow \mathbb{Z}Sing _*(U _0).
\end{equation}
The homology of the total complex of the bicomplex (\ref{***}) gives the Borel-Moore homology \cite{GS}. The weight filtration for $H _*^{BM}(U^{an})$ is the increasing filtration
 $$...\subseteq W _tH _n^{BM}(U^{an})\subseteq W _{t+1}H _n^{BM}(U^{an}) \subseteq ...$$
where
  $$W _tH _n^{BM}(U^{an}):=image(h _n(\mathbb{Z}Sing _*(U _{n+t})\rightarrow ...\rightarrow \mathbb{Z}Sing _*(U _0))\rightarrow H _n^{BM}(U^{an})).$$
It can be proven (\cite{GS}) that 
$$W _tH _n^{BM}(U^{an})=0$$
for any $t<-n$ and
$$W _tH _n^{BM}(U^{an})=H^{BM} _n(U^{an})$$
for any $t\leq d-n$, where by $d$ we denote the dimension of the variety $U$. 

The generalized cycle maps of a quasi-projective variety have the following property:
\begin{proposition}(\cite{FHW})
For any quasi-projective complex variety $U$ the image of the canonical map
$$cyc _{t,n}:L _tH _n(U)\rightarrow H _n^{BM}(U^{an})$$
lies in the part of weight at most $-2t$ of Borel-Moore homology. 
\end{proposition}

 We will recall now the following conjecture due to Bloch and Kato.
\begin{theorem}(Bloch-Kato conjecture)
For any $n\geq 0$ and any field $E$ the norm residue homomorphism 
$$K^M _n(E)/l\rightarrow H^n _{et}(E,\mu^{\tensor q} _m)$$
is an isomorphism.
\end{theorem}
 This conjecture was proven by V. Voevodsky for any $m=2^l$ and for any natural number $l>0$ (this part is also called Milnor conjecture). The general case appears to be proven from the work of M. Rost and V. Voevodsky. A. Suslin and V. Voevodsky \cite{SV2} (see also Geisser-Levine \cite{GL}) proved that the Bloch-Kato conjecture is equivalent to a conjecture due to Beilinson-Lichtenbaum.
\begin{theorem}(Beilinson-Lichtenbaum conjecture)
The map
$$H^n _\mathbb{M}(X, \mathbb{Z}/m(q))\rightarrow H^n _{et}(X,\mu^{\tensor q} _m)$$
is isomorphism for $n\leq q$ and injective for $n\leq q+1$ for any $X$ smooth quasi-projective variety.
\end{theorem}
  A. Suslin proposed the following characterization of morphic cohomology with integral coefficients (see \cite{FHW} and \cite{W}).

\begin {conjecture} (Suslin's conjecture)
The map
$$L^qH^n(X, \mathbb{Z})\rightarrow H^n(X,\mathbb{Z})$$
is isomorphism for $n\leq q$ and injective for $n\leq q+1$ for any $X$ smooth quasi-projective variety.
\end{conjecture}
We notice that the last conjecture contains a conjecture due to E. Friedlander and B. Mazur \cite{FM}.

\begin{conjecture}(Friedlander-Mazur conjecture)
For any complex smooth quasi-projective variety $X$
$$L^qH^n(X)=0$$
for any $n<0$.
\end{conjecture}
 The Friedlander-Lawson duality theorem (\cite{FL1}, \cite{F3}) between morphic cohomology and Lawson homology will be used throughout the paper.
\begin{theorem}(Friedlander-Lawson \cite{FL1}, Friedlander \cite{F3}) For any $X$ quasi-projective smooth complex variety of dimension $d$ 
$$L^sH^n(X)\simeq L _{d-s}H _{2d-n}(X)$$
for any $n\leq 2s$, $n\in \mathbb{Z}$ and $0\leq s\leq 2d$.
\end{theorem}
\section{First results concerning generalized cycle maps}
    We start this section with some applications of the Bloch-Kato conjecture in the context of Lawson homology. The point $b)$ in proposition \ref{p1} is known as Bloch's theorem \cite{B}. In proposition \ref{p6} we analize the torsion of the Borel-Moore homology of a smooth quasi-projective variety.  
\begin{proposition}
\label{p1}
Let $X$ a quasi-projective smooth variety. Then:     

a) Let $n\leq q+1$. Then $ K^{q,n}$ is divisible and $C^{q,n}$ is torsion free. 

b) Suppose $X$ is projective. Then the torsion of $H^n(X)$ is supported in codimension one for any $n>0$.

c) $L^qH^n(X)$ is uniquely divisible for $n<0$ and $L^qH^0(X)$ is torsion free (for any $q\geq 0$). 
\end{proposition}

\begin{proof}
We write  the diagram of universal coefficient sequences for both cohomologies:

$$\begin{CD}
0@>>>L^qH^n(X)\tensor {\mathbb{Z}/m} @>>> L^qH^n(X,\mathbb{Z}/m) @>>>  _mL^qH^{n+1}(X)  @>>> 0\\
 @.    @VVV  @VVV  @VVV\\
0@>>>  H^n(X)\tensor \mathbb{Z}/m  @>>> H^n(X,\mathbb{Z}/m) @>>>  _mH^{n+1}(X)@>>> 0.\\
\end{CD}$$

We recalled in the second section that Bloch-Kato conjecture implies that the map
$$H^n _\mathbb{M} (X, \mathbb{Z}/m(q))\rightarrow H^n _{et}(X,\mu^{\tensor q} _m)$$
is isomorphism for $n\leq q$ and injective for $n\leq q+1$ for any $X$ smooth quasi-projective variety. The above map factors through the cycle map from morphic cohomology to the singular cohomology \cite{FL}. In \cite{SV1} it is proved that 
$$H^n _\mathbb{M} (X, \mathbb{Z}/m(q))\simeq L^qH^n(X,\mathbb{Z}/m)$$
for any allowed $n,q$ and any complex projective variety $X$.

In conclusion the middle vertical map from the above diagram is injective for $n\leq q+1$ and isomorphism for $n\leq q$. Using  the snake lemma we conclude that we have the following exact sequence:
$$
 0\rightarrow L^qH^n(X)\tensor \mathbb{Z}/m\rightarrow H^n(X)\tensor \mathbb{Z}/m \rightarrow Ker( _mL^qH^{n+1}(X)\rightarrow _mH^{n+1}(X))\rightarrow 0
$$
 for any $n\leq q$ and that the map $L^qH^n(X)\tensor \mathbb{Z}/m\rightarrow H^n(X)\tensor \mathbb{Z}/m$ is an injection for $n\leq q+1$. Moreover we conclude that for $n\leq q$ we have 
$$
 _mL^qH^{n+1}(X)\rightarrow _mH^{n+1}(X)
$$
surjective. This means actually that 
$$
_mL^qH^n(X)\rightarrow  _mH^n(X)
$$
is surjective for any $n\leq q+1$ and any $m>1$. It implies that 
\begin{equation}
\label{****}
torsion(Im(cyc^{q,n}))=torsion(H^n(X))
\end{equation}
for any $n\leq q+1$. Because each image of the cycle map is included in a step of the coniveau filtration we have  
$$torsion(N^{n-q}H^{n}(X))=torsion(H^n(X))$$ 
for any $n\leq q+1$. The only case when we conclude something non-trivial from the equality above is when $n=q+1$. In this case for any $0<n\leq d+1$ we have 
$$torsion(N^1H^{n}(X))=torsion(H^n(X))$$
which implies our point b).
 Consider now the composition 
$$L^qH^n(X)\tensor \mathbb{Z}/m\rightarrow Im(cyc^{q,n})\tensor \mathbb{Z}/m\rightarrow H^n(X)\tensor \mathbb{Z}/m.
$$
 The first map is still surjective because $-\tensor \mathbb{Z}/m$ is a right exact functor. For $n\leq q+1$ the composition is injective. This implies that 
 $$L^qH^n(X)\tensor \mathbb{Z}/m\simeq Im(cyc^{q,n})\tensor \mathbb{Z}/m$$
and that
$$Im(cyc^{q,n})\tensor \mathbb{Z}/m\hookrightarrow H^n(X)\tensor \mathbb{Z}/m$$
for any  $n\leq q+1$ and $m>1$.
Consider now the following short exact sequence
$$ 0\rightarrow K^{q,n}\rightarrow L^qH^n(X)\rightarrow Im(cyc^{q,n})\rightarrow 0.$$
Tensoring with $\mathbb{Z}/m$ we obtain the following exact sequence:
$$ 0\rightarrow  _mK^{q,n}\rightarrow  _mL^qH^n(X)\stackrel{a _1}{\rightarrow}   _mIm(cyc^{q,n})\rightarrow K^{q,n}\tensor \mathbb{Z}/m$$
$$\rightarrow L^qH^n(X)\tensor \mathbb{Z}/m\stackrel{a _2}{\rightarrow} Im(cyc^{q,n})\tensor \mathbb{Z}/m\rightarrow 0.$$
For $n\leq q+1$ the map $a _2$ is an isomorphism and the map $a _1$ is a surjection. From exactness of the sequence we get 
$$K^{q,n}\tensor \mathbb{Z}/m=0$$
for any $n\leq q+1$ and $m>1$. This implies that $K^{q,n}$ is divisible  for $n\leq q+1$. 
 
Consider now the following exact sequence
$$0\rightarrow Im(cyc^{q,n})\rightarrow H^n(X)\rightarrow C^{q,n}\rightarrow 0.$$
Tensoring with $\mathbb{Z}/m$ we obtain the following long exact sequence:
$$0\rightarrow _mIm(cyc^{q,n})\stackrel{a _3}{\rightarrow} _mH^n(X)\rightarrow _mC^{q,n}\rightarrow Im(cyc^{q,n})\tensor \mathbb{Z}/m\stackrel{a _4}{\rightarrow} H^n(X)\tensor \mathbb{Z}/m\rightarrow C^{q,n}\tensor \mathbb{Z}/m\rightarrow 0.$$
For $n\leq q+1$ the map $a _3$ is bijective and the map $a _4$ is injective. From the exactness of the sequence we get 
$$ _mC^{q,n}=0$$
for any $n\leq q+1$ and any $m>1$. This implies that $C^{q,n}$ is torsion free for any $n\leq q+1$.

Suppose now that $n<0$. Because $0\leq q \leq d=dim(X)$ we have $n<q$. We have the following short exact sequence
$$
 0\rightarrow L^qH^n(X)\tensor \mathbb{Z}/m\rightarrow L^qH^n(X,\mathbb{Z}/m)\rightarrow _mL^qH^{n+1}(X)\rightarrow 0.$$
Because $L^qH^n(X,\mathbb{Z}/m)=0$ for any $n<0$ and for any $m>1$ we conclude that $L^qH^n(X)\tensor \mathbb{Z}/m=0$ for any $n<0$, $m>1$ (i.e. $L^qH^n(X)$ is divisible for $n<0$) and that $_mL^qH^{n+1}(X)=0$ (i.e. $L^qH^n(X)$ is torsion free for any $n\leq 0$). 

\end{proof}
\begin{corollary}
\label{c1}
Let $n\leq q+1$. Then:

a) If  there is a natural nonzero number $M$ such that $MK^{q,n}=0$ then $cyc^{q,n}$ is injective.

b) Suppose $cyc^{q,n}\tensor \mathbb{Q}$ is surjective. Then $cyc^{q,n}$ is surjective.
\end{corollary}
\begin{proof}
a) From the above proposition \ref{p1} we get $K^{q,n}$ divisible. This implies that for any $x\in K^{q,n} $ there is an element $y\in K^{q,n}$ such that $x=My=0$.

b) From $cyc^{q,n}\tensor \mathbb{Q}$ surjective it follows that $C^{q,n}\tensor \mathbb{Q}=0$. From proposition \ref{p1} we know that $C^{q,n}$ is torsion free. This implies that $C^{q,n}=0$.
\end{proof}
Point b) in proposition  \ref{p1} has the following formulation in the quasi-projective case:
\begin{proposition}
\label{p6}
Let $X$ a smooth quasi-projective variety of dimension $d$ and let $n\geq d$. Fix $s=n-d+1$. Then 
$$ torsion(W _{-2s}H _n(X))=torsion(H^{BM} _n(X))$$
where $W _{-2s}H _n(X)$ is a step in the weight filtration of the Borel-Moore homology $H^{BM} _n(X)$ \cite{FHW}.
\end{proposition}
\begin{proof}
From (\ref{****}) we have that the groups $Im(cyc^{q,q+1})$ and $H^{q+1}(X)$ have the same torsion for any $q$ with $0\leq q\leq d$. Using the Friedlander-Lawson duality theorem \cite{F3} we have that 
$$tors(Im(L _{n-d+1}H _n(X)\rightarrow H _n(X)))=tors(H _n(X))$$ 
for any $n\geq d$.

 We recalled in the second section that the cycle map from Lawson homology to the Borel-Moore homology of a smooth quasi-projective variety factors through steps in the weight filtration \cite{FHW}, i.e.
  $$L _sH _n(X)\rightarrow W _{-2s}H _n(X)\hookrightarrow H _n(X)$$
for any $0\leq s\leq d$ and $n\geq 2s$. This implies the statement of the theorem.
\end{proof}
 
 The above discussion gives us the following reformulation of the Friedlander-Mazur conjecture.
\begin{proposition}
Let $X$ be a smooth quasi-projective variety. Then the Friedlander-Mazur conjecture is valid for $X$ if and only if $L^qH^n(X) _\mathbb{Q}\simeq 0$ for any $n<0$.
\end{proposition}
\begin{proof}
The point c) in the proposition \ref{p1} shows that these groups are torsion free. 
\end{proof}

\begin{remark}(cohomological Brauer group)

Let $X$ a smooth projective variety of dimension $d$ with $$H^2 _{Zar}(X,O _X)=H^1 _{Zar}(X,O _X)=0.$$ We know (see for example \cite{Da}) that, in these conditions, the cohomological Brauer group of $X$ has the following characterization
$$Br(X)\simeq tors(H^3(X)).$$
Suslin's conjecture predicts that the cycle map 
$$L^qH^n(X)\rightarrow H^n(X)$$
is an isomorphism for any $n\leq q$ and a monomorphism for $n=q+1$.
Let's assume Suslin's conjecture. We obtain that
$$torsion(L^3H^3(X))\simeq torsion(L^4H^3(X))\simeq ...\simeq torsion(H^3(X)).$$
But (\ref{****}) shows that 
$$torsion(Im(cyc^{2,3}))=torsion(H^3(X))$$
and because Suslin's conjecture gives us that the cycle map $cyc^{2,3}$ is injective we obtain that 
$$torsion(L^2H^3(X))\simeq torsion(L^3H^3(X))\simeq...\simeq torsion(H^3(X))$$
giving a characterization of the cohomological Brauer group of $X$ by means of morphic cohomology. We will show in sections five, six and seven that Suslin's conjecture can be verified for certain projective varieties. 

 A natural question to ask is whether $tors(L^2H^3)$ is a birational invariant in general as Suslin's conjecture predicts. We will prove below that this is indeed the case. We will use a blow-up formula for Lawson homology proved by Hu \cite{Hu} and the fact that birational maps between projective smooth varieties factor as a composition of blow-ups with centers of codimension greater that two \cite{Ab}. 
\end{remark}
\begin{proposition}
$tors(L^2H^3)$ is a birational invariant. 

That is for any $X$, $X^\prime$ birational equivalent smooth projective varieties we have:
$$tors(L^2H^3(X))\simeq tors(L^2H^3(X^\prime)).$$
\end{proposition}
\begin{proof}
Let $X _Y\rightarrow X$ a blow up of a smooth center $Y$ of codimension greater or equal with 2. From \cite{Hu} we know that
$$L _{n-2}H _{2n-3}(X _Y)=\oplus _{1\leq j\leq r-1}L _{n-2-j}H _{2n-3-2j}(Y)\oplus L _{n-2}H _{2n-3}(X)$$
where by $r$ we denote the codimension of $Y$ in $X$.

 It suffices to show that $tors (L _{n-2-j}H _{2n-3-2j}(Y))=0$ for any $1\leq j\leq r-1$. We notice that 
$$dimY-1\leq n-2-j\leq n-3$$
and 
 $$2dim(Y)-1\leq 2n-3-2j\leq 2n-5.$$
 If $n-2-j\geq dim(Y)$ it is obvious that 
$$tors(L _{n-2-j}H _{2n-3-2j}(Y))=0.$$
 We also have 
$$tors(L _{dimY-1}H _{2dimY-1}(Y))\simeq tors(H _{2dimY-1}(Y)).$$
 From the universal coefficient sequence we obtain 
$$tors(H _{2dimY-1}(Y))=tors(H _0(Y))=0.$$
We can conclude now that 
 $$tors(L^2H^3(X))\simeq tors(L^2H^3(X^\prime))$$  
for any $X,X^\prime$ birational equivalent smooth projective varieties.
\end{proof}
\section{Cycle action on morphic cohomology}
 Let $\alpha$ be a dimension $d=dim(X)$ irreducible algebraic cycle in $ X\times X$ with the support contained in $V\times W$, where $V\subset X$ and $W\subset X$ are irreducible subvarieties. Denote $v=dim(V)$ and $w=dim(W)$. Consider the resolutions of singularities $i:V^*\rightarrow V$ and $j:W^*\rightarrow W$. Denote $\alpha^*\in CH _d(V^*\times W^*)$ an element such that 
$$(i\times j) _*\alpha^*=\alpha.$$
We remark that we can always find such a cycle up to a moving in the rational equivalence class of $\alpha$ \cite {EL}. 
The cycle $\alpha$ gives the following action
$$\alpha _*:L^mH^l(X)\rightarrow L _{d-m}H _{2d-l}(X)$$
with $d=dim(X)$ and 
$$\alpha _*(x)=pr _{2*}(pr _1^*(x)\cap\alpha)$$
The above map depends only on the algebraic equivalence class of $\alpha$. A similar action in the context of Lawson homology was also considered by C. Peters in \cite{P} (see also \cite{FL}). 

 The above map decomposes in the following way
$$\alpha _*(x)=j _*pr _{W^**}(pr _{V^*}^*i^*(x)\cap \alpha^*)$$
by using projection formula (see \cite{BO} for a proof of the projection formula) in the morphic cohomology setting. If we consider for example the action of the diagonal cycle of
$X\times X$ we will obtain the Friedlander-Lawson duality isomorphism between
morphic cohomology and  Lawson homology.  

 The above action commutes with the similarly defined action of the algebraic cycle
$\alpha$ on the singular cohomology. This can be seen from the fact
that the cycle maps from morphic cohomology (resp. Lawson homology) to
singular cohomology (resp. singular homology) are natural \cite{BO} and
commute with cap product with an algebraic cycle \cite{BO}.  

The above discussion is summarized in the following sequence of the
commutative diagrams (the horizontal maps are given by the
decomposition of the actions of the algebraic cycle $\alpha$ and the vertical
maps are given by the cycle maps)  

$$\begin{CD}
 L^{d-m}H^{2d-l}(X)@>i^*>> L^{d-m}H^{2d-l}(V^*) @>pr^* _{V^*}>> L^{d-m}H^{2d-l}(V^*\times W^*) @>\cap \alpha^*>> @.\\
 @Vc _1VV    @Vc _2VV  @VVV  \\
 H^{2d-l}(X)@>i^*>> H^{2d-l}(V^*)  @>pr^* _{V^*}>>H^{2d-l}(V^*\times W^*)   @>\cap \alpha^*>>  @.\\
\end{CD}$$
\begin{equation}
\begin{CD}
\label{dia1}
 @>\cap \alpha^*>> L _mH _l(V^*\times W^*) @>(pr _{W^*}) _*>> L _mH _l(W^*) @>j _*>> L _mH _l(X)\\
   @.         @VVV            @Vc _3VV         @Vc _4VV\\
   @>\cap \alpha^*>>   H _l(V^*\times W^*) @>(pr _{W^*}) _*>> H _l(W^*) @>j _*>> H _l(X)\\
\end{CD}
\end{equation}
for any $0\leq m\leq d$ and any $l\geq 2m$.

The map $c _2$ is an isomorphism for any $m\leq d-v$, where $v=dim(V)$. To see this we divide it in several cases depending on the value of $2d-l$. If $2d-l>2v\geq 0$ then   $$L^{d-m}H^{2d-l}(V^*)= H^{2d-l}(V^*)=0$$
by \cite{FL}. If $0\leq 2d-l\leq 2v$ then we can consider the following morphic cohomology group $L^vH^{2d-l}(V^*)$ which is isomorphic with $H^{2d-l}(V^*)$ from the Poincare duality and the Dold-Thom theorem (see \cite{BO}). At the same time the composition of s-maps$$ L^vH^{2d-l}(V^*)\rightarrow L^{d-m}H^{2d-l}(V^*)$$
is an isomorphism in this range and commutes with the cycle maps \cite{FL1}. This implies that 
$$L^{d-m}H^{2d-l}(V^*)\simeq H^{2d-l}(V^*).$$
Consider now $2d-l<0$. Then by the Friedlander-Lawson duality theorem and the isomorphism of s-maps in this range we obtain
$$0=L^vH^{2d-l}(V^*)=L^{d-m}H^{2d-l}(V^*).$$
The map $c _3$ is an isomorphism for any $m\geq w$. In the case $m=w$ we have $$L _mH _{2m}(W^*)\simeq H _{2m}(W^*)$$ because $W^*$ is irreducible. For $m>w$ we obviously have $L _mH _l(W^*)=H _l(W^*)=0$ since $l\geq 2m>2w$.

The above discussion proves the following proposition:
\begin{proposition}
\label{l1}
 Let $\alpha$ be an irreducible algebraic cycle in $CH^d(X\times X)$ with the support contained in $V\times W$, where $V\subset X$ and $W\subset X$ are irreducible subvarieties of dimension $v$, respectively $w$.
The action of the cycle $\alpha$ on the kernel and the cokernel of the map $c _1$ is zero for $m\geq w=dim(W)$ or for $m\leq d-v=codim(V)$.
\end{proposition}
\begin{proof}
From the above discussion we conclude that if $m\geq w$ then $c _3$ is an isomorphism and that if $m\leq d-v$ then $c _2$ is an isomorphism. These imply the conclusion of our proposition.
\end{proof}
\begin{corollary}
\label{l2}
Suppose $\alpha$ as in proposition \ref{l1} and suppose $dim(X)=dim(V)+dim(W)$. Then the action of the cycle $\alpha$ on the kernel and on the cokernel of the map $c _1$ is zero for any $0\leq m\leq d$.
\end{corollary}
\begin{proof}
Direct consequence of proposition \ref{l1}.
\end{proof}
\begin{remark}
\label{r9}
We remark that to study the action of a cycle $\alpha =\sum n _i\alpha _i\in CH^d(X\times X)$ with $supp(\alpha) \subset V\times W$ it is enough to study the action of each irreducible cycle $\alpha _i$. It is obvious that 
$$ supp(\alpha _i)\subset supp(\alpha) \subset V\times W$$
and that because $supp(\alpha _i)$ is irreducible there are $V _i\subset V, W _i\subset W$ irreducible components such that 
$$ supp(\alpha _i)\subset V _i\times W _i$$
\end{remark}
By using the Friedlander-Lawson duality theorem \cite{FL1} we will identify the cycle map $c_1$ with the cycle map $L _mH _l(X)\rightarrow H _l(X)$.

This will identify the action of the diagonal cycle with the identity map.

\begin{convention}
\label{co1}
From now on by ``the action of $\alpha$ is zero for $m$ in some certain range'' we will understand that the action of the cycle $\alpha$ on the kernel and cokernel of the cycle map $L _mH _*\rightarrow H _*$ is zero for $m$ in the respective range. 
\end{convention}
\section{Comparing Lawson homology with singular homology}
In this section we study the cycle maps $cyc^{q,n}:L^qH^n(X)\rightarrow H^n(X)$ for $X$ smooth projective complex variety with the property that its zero cycles are supported on a proper subvariety. We prove that these cycle maps behave nicely for threefolds and fourfolds with this property (being most of the time injective or bijective). We expect that there are cycle maps $cyc^{q,n}$ totally non-trivial for varieties $X$ of large dimension with zero cycles supported on a subvariety. As a support for our expectation is a theorem of A. Albano, and A. Collino \cite{AC} proving that for a generic smooth cubic hypersurface $X\subset \mathbb{P}^8$ the Griffiths group $Griff^4(X)\tensor \mathbb{Q}$ is infinitely generated.

We start the section by recalling a result of E. Friedlander. He proved \cite{F} that for any smooth connected complex projective variety $X$ of dimension $d$ we have
$$L _{d-1}H _{2d-2}(X)\hookrightarrow H _{2d-2}(X),$$
$$L _{d-1}H _{2d-1}(X)\simeq H _{2d-1}(X),$$
$$L _{d-1}H _{2d}(X)\simeq H _{2d}(X)\simeq \mathbb{Z}$$
and that 
$$L _{d-1}H _k(X)=0$$
for any $k>2d$.

We recall that Bloch and Srinavas \cite{BS} proved that if a smooth projective variety $X$ has its zero cycles supported on a not necessary ireducible subvariety $V$, i.e. $CH _0(X\setminus V)=0$, then the diagonal cycle decomposes as
$$N\Delta=\alpha +\beta$$
for some natural nonzero number $N$ and some cycles $\alpha,\beta\in CH^d(X\times X)$ with the support of $\alpha$ included in $V\times X$ and the support of $\beta$ included in  $X\times D$, where $D$ is a divisor of $X$. We will also use the transpose of this decomposition, i.e
  $$N\Delta=\alpha ^t +\beta ^t$$
where $\alpha ^t ,\beta ^t \in CH^d(X\times X)$ are supported on $X\times V$, respectively $D\times X$.

The following theorem computes $K^{sst}$ for ``degenerate'' threefolds.

\begin{theorem}
\label{t1}
Let $X$ be a smooth projective complex threefold such that there is a proper subvariety $V\subset X$ with $CH _0(X \setminus V)=0$. Then:
$$K^{sst} _i(X)\simeq ku^{-i}(X^{an}), i\geq 1,$$
$$K^{sst} _0(X)\hookrightarrow ku^0(X^{an}).$$
Moreover if $X$ is a rationally connected threefold then
$$K^{sst} _i(X)\simeq ku^{-i}(X^{an})$$
for any $i\geq 0$.
\end{theorem}

Some examples which fulfill the conditions above are: rationally connected threefolds (e.g. smooth Fano threefolds \cite{KMM}), Kummer threefolds \cite{BS}, certain quotient varieties such $$(X\times E)/(\mathbb{Z}/2)$$ with $X$ a K3 covering of an Enriques surface  and $E$ an elliptic curve \cite{BS}.

The above result on rationally connected threefolds generalizes the same result on rational threefolds proved in \cite{FHW} with other tools.
\begin{proof}
The proof of the above theorem is based on the spectral sequence relating morphic cohomology and semi-topological K-theory \cite{FHW} and on the following two propositions which compute the Lawson homology groups of a threefold $X$ with zero cycles supported on a subvariety.
\begin{proposition}
\label{p2}
Let $X$ be a smooth projective complex threefold such that there is a proper subvariety $V\subset X$ with $CH _0(X \setminus V)=0$ and $dim(V)\leq 1$. Then:

a) $L _1H _2(X)\hookrightarrow H _2(X)$ is injective and a rational isomorphism.

b) $L _1H _3(X)\simeq H _3(X)$.

c) $L _2H _4(X)\simeq L _1H _4(X)\simeq H _4(X)$.

d) $L _2H _5(X)\simeq L _1H _5(X)\simeq H _5(X)$.

e) $L _3H _6(X)\simeq L _2H _6(X)\simeq L _1H _6(X)\simeq H _6(X).$

f) $L _kH _n(X)=0$ for any $n\geq 7$ and any $k\geq 0$.

In particular any such threefold fulfills Suslin's conjecture.

Moreover if $X$ is rationally connected threefold then 
$$L _*H _*(X)=H _*(X)$$
for all possible indices.

\end{proposition}
\begin{proof} We proof the case $dim(V)=1$, the other case being similar.

Consider the above decomposition $$N\Delta=\alpha ^t +\beta ^t$$ with $\alpha ^t$ supported on $X\times V$ and $\beta ^t$ supported on $D\times X$, with $D$ a divisor in $X$. Remark \ref{r9} shows that it is enough to consider the case when $V$ and $D$ are irreducible. We recall that we denote by $D^*$, respectively $V^*$, the resolution of singularities of $D$, respectively $V$. 

Proposition \ref{l1} gives us that the action of $\beta^t$ is zero on Ker$(L _mH _*(X)\rightarrow H _*(X))$ for $m\leq codim(D)=1$ and the action of $\alpha ^t$ on the same kernel is zero for $m\geq v=dim(V)=1$. This implies that for $m=1$ we have
$$N(Ker(L _1H _k(X)\rightarrow H _k(X)))=0$$
for $2\leq k\leq 6$ and that 
$$NL _1H _k(X)=0$$
for $k\geq 7$. Proposition \ref{p1} implies that $L _1H _k(X)=0$ for any $k\geq 7$ and Corollary \ref{c1} implies that the cycle map $L _1H _k(X)\rightarrow H _k(X)$ is injective for any $3\leq k\leq 6$.

Let $x\in H _k(X)\simeq H^{6-k}(X)$ with $3\leq k\leq 6$. Then $$\beta ^t _*x\in Im(L _1H _k(X)\rightarrow H _k(X))$$ because $L^2H^{6-k}(D^*)\simeq H^{6-k}(D^*)$ (see diagram (\ref{dia1})). 

We have $\alpha ^t _*x=0$ because the action of $\alpha ^t$ on $H _k(X)$ factors through $H _k(V^*)$, $dim(V)=1$ and $ k\geq 3$. This implies that 
$$Nx=\beta ^t _*x\in Im(L _1H _k(X)\rightarrow H _k(X))$$
for any $3\leq k\leq 6$. Because $x\in H _k(X)$ is arbitrary we conclude that the rational cycle map   
$$L _1H _k(X)\tensor \mathbb{Q}\rightarrow H _k(X)\tensor \mathbb{Q}$$
is surjective for any $3\leq k\leq 6$. Corollary \ref{c1} shows that the cycle map $L _1H _k(X)\rightarrow H _k(X)$ is surjective for any $3\leq k\leq 6$.

For $k=2$ we use a result of Bloch-Srinavas. They prove that for varieties as in our hypothesis algebraic equivalence and homological equivalence coincides for codimension 2 cycles \cite{BS}. This means that the cycle map
$$L _1H _2(X)\rightarrow H _2(X)$$
is injective.

Consider now the decomposition  $$N\Delta=\alpha +\beta $$ with $\alpha$ supported on $V\times X$ and $\beta $ supported on $X\times D$. Proposition \ref{l1} gives that the action of $\alpha$ is zero on Ker$(L _mH _*(X)\rightarrow H _*(X))$ for $m\leq d-v=codim(V)=2$ and that the action of $\beta$ is zero on Ker$(L _mH _*(X)\rightarrow H _*(X))$ for $m\geq dim(D)=2$. This implies that for $m=2$ we have
$$N(Ker(L _2H _*(X)\rightarrow H _*(X)))=0$$
for $4\leq k\leq 6$ and that 
$$NL _2H _k(X)=0$$
for $k\geq 7$. Proposition \ref{p1} implies that $L _2H _k(X)=0$ for any $k\geq 7$ and Corollary \ref{c1} implies that the cycle map $L _2H _k(X)\rightarrow H _k(X)$ is injective for any $4\leq k\leq 6$. 

Let $x\in H _k(X)$ with $4\leq k\leq 6$. Then $$\alpha _*(x)\in Im(L _2H _k(X)\rightarrow H _k(X)) $$ because $L^1H^{6-k}(V^*)\simeq H^{6-k}(V^*)$ (see diagram (\ref{dia1})).
The action of $\beta$ on $H _k(X)$ is zero for $5\leq k\leq 6$ because this action factors through $H _k(D^*)$ and $dim(D^*)=2$. If $k=4$ then $$\beta _*(x)\in Im(L _2H _4(X)\rightarrow H _4(X))$$ because $L _2H _4(D^*)\simeq H _4(D^*)$ (see diagram (\ref{dia1})).  Because 
$$Nx=\alpha _*(x)+\beta _*(x)$$
we conclude that for $4\leq k\leq 6$ the cycle maps
$$L _2H _k(X)\tensor \mathbb{Q}\rightarrow H _k(X)\tensor \mathbb{Q}$$
are surjective. Applying Corollary \ref{c1} we conclude that these surjections are with integer coefficients.

Let $X$ be a rationally connected threefold. Then $H^4(X,\mathbb{C})=H^{2,2}(X)$ because $h^{1,3}=h^{3,1}=h^{2,0}=0$. This implies that 
$$ H _2(X)\simeq H^4(X)\simeq H^{2,2}(X,\mathbb{Z})$$
where we denoted $H^{2,2}(X,\mathbb{Z}):=\{\eta\in H^4(X)$ such that $coef _*(\eta)\in H^{2,2}(X)$ with $coef _*: H^4(X,\mathbb Z)\rightarrow H^4(X,\mathbb C)$ being the coefficient map $\}$.

C. Voisin proved the following theorem:
\begin{theorem}(Voisin \cite{V})

The Hodge conjecture with integral coefficients is valid for any smooth uniruled threefold.

\end{theorem}

This theorem implies that the cycle map $L _1H _2(X)\rightarrow H _2(X)$ is surjective for $X$ any smooth rationally connected threefold.
\end{proof}

\begin{proposition}
\label{p3}
Let $X$ be a smooth projective complex threefold such that there is a proper subvariety $V\subset X$ with $CH _0(X \setminus V)=0$ and $dim(V)=2$. Then:

a) $L _1H _2(X)\hookrightarrow H _2(X)$.

b) $L _1H _3(X)\simeq H _3(X)$.

c) $L _2H _4(X)\hookrightarrow L _1H _4(X)\simeq H _4(X).$

d) $L _2H _5(X)\simeq L _1H _5(X)\simeq H _5(X)$.

e) $L _3H _6(X)\simeq L _2H _6(X)\simeq L _1H _6(X)\simeq H _6(X).$

f) $L _kH _n(X)=0$ for any $n\geq 7$ and any $k\geq 0$.

In particular any such threefold fulfills Suslin's conjecture.
\end{proposition}

\begin{proof}
Consider the decomposition  $$N\Delta=\alpha +\beta $$ with $\alpha$ and $\beta$ being supported on $V\times X$, respectively $X\times D$. It is enough to consider the case when $V$ and $D$ are irreducible (see remark \ref{r9}). The action of $\alpha$ is zero for $m\leq d-v=codim(V)=1$ and the action of $\beta$ is zero for  $m\geq dim(D)=2$ (see Convention \ref{co1}). 

Suppose now $m=1$. Because $D^*$ is a surface we have for $l\geq 3$
$$L _1H _l(D^*)\simeq H _l(D^*).$$
This implies that the action of $\beta$ is zero for $m=1$ and $l\geq 3$ (see Diagram \ref{dia1}). Because we already know that the action of  $\alpha$ is zero for $m=1$ it implies that 
$$L _1H _l(X) _\mathbb{Q}\simeq H _l(X) _\mathbb{Q}$$
for any $l\geq 3$. Applying corollary \ref{c1} we obtain 
$$L _1H _l(X)\simeq H _l(X)$$
for any $l\geq 3$.

Suppose now $m=2$. Because $V^*$ is a surface we get
$$L^1H^{6-l}(V^*)\simeq H^{6-l}(V^*)$$
for any $l\geq 5$ which means that for these indexes the action of $\alpha$ is zero. As we know that for $m=2$ the action of $\beta$ is zero we can conclude as above that
$$L _2H _l(X)\simeq H _l(X)$$
for any $l\geq 5$.

The injectivity with integer coefficients in a) and c) follows from the theorem of Bloch-Srinavas \cite{BS} used in the proposition \ref{p2} and from the fact that for divisors algebraic equivalence coincides with homological equivalence.
\end{proof}
Using now the spectral sequence argument from \cite{FHW}, theorem 6.1 and our propositions \ref{p2} and \ref{p3} we can conclude our theorem \ref{t1}.

\end{proof}

\begin{remark}
\label{r1}
Bloch-Srinavas \cite{BS} proved that a Kummer threefold has its zero cycles supported on a subvariety of dimension two. These Kummer threefolds show that the injectivity in the point a) of proposition \ref{p3} is the best  we can get because a Kummer threefold has $h^{2,0}\neq 0$.
\end{remark}
\begin{remark}
For $X$ and $V$ as in the theorem \ref{t1} and $dim(V)\leq 1$ the injection $$K^{sst} _0(X)\hookrightarrow ku^0(X^{an})$$ is moreover a rational isomorphism. For $dim(V)=2$ this rational isomorphism is not longer true as we can see from the example given in remark \ref{r1}.
\end{remark}
\begin{convention}
 Remark \ref{r9} shows that it is enough to study the action of a irreducible cycle. In the rest of the paper, without reducing the generality, we will understand that a decomposition of the form 
$$N\Delta=\alpha +\beta$$
with $\alpha$ supported on $V\times X$ and $\beta$ supported on $X\times D$ has $V$ and $D$ irreducible varieties.
\end{convention}
The next theorem computes $K^{sst}$ for certain ``degenerate'' fourfolds.
\begin{theorem}
\label{t2}
Let $X$ be a smooth projective fourfold such that there is a proper subvariety $V\subset X$ of $dim(V)\leq 2$ with $CH _0(X \setminus V)=0$. Then: 
$$K^{sst} _i(X)\simeq ku^{-i}(X^{an}),i\geq 3,$$
$$K^{sst} _2(X)\hookrightarrow ku^{-2}(X^{an}),$$
$$K^{sst} _i(X) _\mathbb{Q}\simeq ku^{-i}(X^{an}) _\mathbb{Q},i=1,2,$$
$$K^{sst} _0(X) _\mathbb{Q}\hookrightarrow ku^0(X^{an}) _\mathbb{Q}.$$
\end{theorem}
Some examples of varieties which fulfill the conditions of the theorem are: rationally connected fourfolds, certain quotient varieties as in \cite{BS} etc.
\begin{proof}
The proof is similar to the proof of the theorem about $K^{sst}$ for degenerate threefolds. It is a corollary of the spectral sequence relating morphic cohomology groups and $K^{sst}$ and of the computation of some Lawson groups. 
\begin{proposition}
\label{p4}
Let $X$ be a smooth projective fourfold such that there is a proper subvariety $V\subset X$ of $dim(V)\leq 1$ with $CH _0(X \setminus V)=0$. Then:

a)$L _1H _2(X) _{\mathbb{Q}}{\simeq}H _2(X) _{\mathbb{Q}}.$

b)$L _1H _3(X) _{\mathbb{Q}}{\simeq}H _2(X) _{\mathbb{Q}}.$

c)$L _2H _4(X)\hookrightarrow L _1H _4(X)\simeq H _4(X).$

d)$L _2H _5(X)\simeq L _1H _5(X)\simeq H _5(X).$

e)$L _3H _6(X)\simeq L _2H _6(X)\simeq L _1H _6(X)\simeq H _6(X).$

f)$L _3H _7(X)\simeq L _2H _7(X)\simeq L _1H _7(X)\simeq H _7(X).$

g)$L _4H _8(X)\simeq L _3H _8(X)\simeq L _2H _8(X)\simeq L _1H _8(X)\simeq H _8(X).$

h)$L _kH _n(X)=0$ for any $n\geq 9$ and any $k\geq 0$.

In particular any such fourfold fulfills Suslin's conjecture.
\end{proposition}
\begin{proof}
Consider the decomposition  $$N\Delta=\alpha +\beta $$ with $\alpha$ supported on $V\times X$ and $\beta$ supported on $X\times D$. The action of $\alpha$ is zero for $m\leq 3$ and the action of $\beta$ is zero for $m\geq 3$ (see Convention \ref{co1}). This implies that 
$$L _3H _*(X)\tensor \mathbb{Q}\simeq H _*(X)\tensor \mathbb{Q}$$
and because of Corollary \ref{c1} we obtain
$$L _3H _*(X)\simeq H _*(X).$$
Because $D^*$ is a smooth threefold we have that the cycle map
$$L _2H _l(D^*)\rightarrow H _l(D^*)$$
is an isomorphism for $l\geq 5$ and a monomorphism for $l=4$. This implies that the action of $\beta$ on the kernel and cokernel of the cycle map
$$L _2H _l(X)\rightarrow H _l(X)$$
is zero for $l\geq 5$. Because we already know that the action of $\alpha$ is zero for $m=2$ we conclude using Corollary \ref{c1} that 
$$L _2H _l(X)\simeq H _l(X)$$
for any $l\geq 5$. The injection from the point c) comes from the fact that for such varieties algebraic equivalence and homological equivalence coincide on codimension 2 cycles \cite{BS}.

Consider now the decomposition  $$N\Delta=\alpha^t +\beta^t $$ with $\alpha^t$ supported on $X\times V$ and $\beta$ supported on $D\times X$. The action of $\alpha^t$ is zero for $m\geq 1$ and the action of $\beta^t$ is zero for $m\leq 1$ (see Convention \ref{co1}). This implies that 
$$L _1H _*(X)\tensor \mathbb{Q}\simeq H _*(X)\tensor \mathbb{Q}$$
and from Corollary \ref{c1} we obtain
$$L _1H _l(X)\simeq H _l(X)$$
for any $l\geq 4$.
  
\end{proof}
\begin{proposition}
\label{p5}
Let $X$ be a smooth projective fourfold such that there is a proper subvariety $V\subset X$ of $dim(V)=2$ with $CH _0(X \setminus V)=0$. Then:

a)$L _1H _2(X) _{\mathbb{Q}}{\hookrightarrow}H _2(X) _{\mathbb{Q}}.$

b)$L _1H _3(X) _{\mathbb{Q}}{\simeq}H _2(X) _{\mathbb{Q}}.$

c)$L _2H _4(X)\hookrightarrow L _1H _4(X)\simeq H _4(X).$

d)$L _2H _5(X)\simeq L _1H _5(X)\simeq H _5(X).$

e)$L _3H _6(X)\hookrightarrow L _2H _6(X)\simeq L _1H _6(X)\simeq H _6(X).$

f)$L _3H _7(X)\simeq L _2H _7(X)\simeq L _1H _7(X)\simeq H _7(X).$

g)$L _4H _8(X)\simeq L _3H _8(X)\simeq L _2H _8(X)\simeq L _1H _8(X)\simeq H _8(X).$

h)$L _kH _n(X)=0$ for any $n\geq 9$ and any $k\geq 0$.

In particular any such fourfold fulfills Suslin's conjecture.
\end{proposition}
\begin{proof}
Consider the decomposition  $$N\Delta=\alpha +\beta $$ with $\alpha$ supported on $V\times X$ and $\beta$ supported on $X\times D$. The action of $\alpha$ is zero for $m\leq 2$ and the action of $\beta$ is zero for $m\geq 3$  (see Convention \ref{co1}). 

Because $D^*$ is a smooth threefold we have that the cycle map
$$L _2H _l(D^*)\rightarrow H _l(D^*)$$
is an isomorphism for $l\geq 5$ and a monomorphism for $l=4$. This implies that the action of $\beta$ on the kernel and cokernel of the cycle map
$$L _2H _l(X)\rightarrow H _l(X)$$
is zero for $l\geq 5$. Because we already know that the action of $\alpha$ is zero for $m=2$ we conclude using corollary \ref{c1} that 
$$L _2H _l(X)\simeq H _l(X)$$
for any $l\geq 5$. The injection from the point c) comes from the fact that for such varieties algebraic equivalence and homological equivalence coincide on codimension 2 cycles \cite{BS}.

Consider the action of $\alpha$ on the kernel and the cokernel of the cycle maps
$$L _3H _l(X)\rightarrow H _l(X).$$
This action factors through $L^1H^{8-l}(V^*)\simeq L _1H _{l-4}(V^*)$ (see diagram \ref{dia1}). Because $V^*$ is a surface we have that the cycle map
$$L _1H _{l-4}(V^*)\simeq H _{l-4}(V^*)$$
for any $l\geq 7$ and injective for $l=6$. Because the action of $\beta$ is zero for $m=3$ we have that
$$L _3H _l(X)\tensor \mathbb{Q}\simeq H _l(X)\tensor \mathbb{Q}$$
for any $l\geq 7$. From Corollary \ref{c1} we conclude that 
$$L _3H _l(X)\simeq H _l(X)$$
for any $l\geq 7$. The injectivity in point e) comes from the fact that on divisors algebraic equivalence coincides with homological equivalence.

Consider now the decomposition  $$N\Delta=\alpha^t +\beta^t $$ with $\alpha^t$ and $\beta^t$ being supported on $X\times V$, respectively $D\times X$. The action of $\alpha^t$ is zero for $m\geq 2$ and the action of $\beta^t$ is zero for $m\leq 1$. Because $V^*$ is a surface, we have that
$$L _1H _l(V^*)\rightarrow H _l(V^*)$$
is an isomorphism for any $l\geq 3$ and a monomorphism for $l=2$. This implies that 
$$L _1H _l(X)\tensor \mathbb{Q}\rightarrow H _l(X)\tensor \mathbb{Q}$$
is an isomorphism for any $l\geq 3$ and a monomorphism for $l=2$. Using Corollary \ref{c1} we can conclude that 
$$L _1H _l(X)\simeq H _l(X)$$
for any $l\geq 4$.

\end{proof}
Using now the spectral sequence argument from \cite{FHW}, theorem 6.1 and our Propositions \ref{p4} and \ref{p5} we can conclude our Theorem \ref{t2}.

\end{proof}

\begin{remark}
As smooth cubic fourfolds show, the injectivity  
$$K^{sst} _0(X) _\mathbb{Q}\hookrightarrow ku^0(X^{an}) _\mathbb{Q}$$
is the best we can obtain.
\end{remark}
The following proposition was previously known in the case of generic cubic hypersurfaces \cite{I} which are known to be rationally connected. 
\begin{proposition}
Let $X$ be a projective smooth variety of dimension $d\geq 3$ and suppose that there is a subvariety $V\subset X$ of $dim(V)\leq 2$ with $CH _0(X\setminus V)=0$. Then
$$N^1H^d(X)=H^d(X).$$
\end{proposition}
\begin{proof}
Without restricting the generality we may suppose $dim(V)=2$. We consider the following decomposition of the diagonal
$$N\Delta=\alpha ^t +\beta ^t$$ with $\alpha ^t$ supported on $X\times V$ and $\beta ^t$ supported on $D\times X$. As before we obtain that the action of $\alpha ^t$ is zero for $m\geq 2$ and the action of $\beta ^t$ is zero for $m\leq 1$ (see Convention \ref{co1}). Denote $V^*$ a desingularization of $V$. Because $V^*$ is a surface we have 
$$L _1H _l(V^*)\tensor \mathbb{Q}\simeq H _l(V^*)\tensor \mathbb{Q}$$
for any $l\geq 3$. This implies that $$L _1H _l(X) _\mathbb{Q}\simeq H _l(X) _\mathbb{Q}$$ 
for any $l\geq 3$. From Corollary \ref{c1} we get 
$$L _1H _d(X)\simeq H _d(X)$$
and because the image of this cycle map is included in $N _{d-1}H _d(X)=N^1H^d(X)$ we obtain our conclusion.
\end{proof}
 
\section{About varieties with small Chow group}
It is proved by Jannsen \cite{Ja} and Laterveer \cite{Lat} that if the the following cycle maps are injective
$$CH _k(X)\tensor \mathbb{Q}\hookrightarrow H _{2k}(X)\tensor \mathbb{Q}$$
for any $0\leq k\leq r$ then we have the following decomposition of the diagonal
\begin{equation}
\label{ec}
N\Delta=\alpha _0+\alpha _1+\alpha _2+...+\alpha _r+\beta
\end{equation}
where $\alpha _i$ are supported on $V _i\times W _{d-i}$ and $\beta$ is supported on $X\times \Gamma^{r+1}$ and $N$ is a nonzero natural number (the lower indices represent the dimension of the subvariety and the upper indices represent the codimension of the subvariety). We denote by $d$ the dimension of the variety $X$.

We say that $X$ has small Chow group of rank $r$ if and only if the first $r$ cycle maps are injective (\cite{EL}, \cite{Lat}).
The following theorem and corollary are extensions of the main results of C. Peters \cite{P}.
\begin{theorem}
\label{t3}

Let $X$ be a smooth projective variety with small Chow group of rank $r$. Then 
there is a natural nonzero number $N$ such that 

a)$NK _{s,*}=0$ for any $s\in \{0,1,..,r+1\}$ (\cite{P}). 

b)$NK^{s,*}=0$ for any $s\in \{0,1,..,r+2\}.$
\end{theorem}
\begin{proof}
Because $X$ has small Chow group of rank $r$, the diagonal cycle decomposes as in \ref{ec}. Because the cycles $\alpha _i$ are supported on $V _i\times W _{d-i}$ with $dim(V _i)+dim(W _{d-i})=dim(X)$ we know from Corollary \ref{l2} that the action of $\alpha _i$ is zero for any $m$ (see Convention \ref{co1}). This implies that
$$N\Delta _*=\beta _*$$
on the kernel and the cokernel of the cycle map
$$L _mH _l(X)\rightarrow H _l(X).$$
But the action of $\beta$ is zero for $m\geq n-r-1$ because it factors through $L _mH _l(\Gamma^{r+1})$ (see Diagram \ref{dia1}). We know \cite{F} that the cycle map
$$L _{n-r-2}H _*(\Gamma^{r+1*})\rightarrow H _*(\Gamma^{r+1*})$$
is injective (where $\Gamma^{r+1*}\rightarrow \Gamma^{r+1}$ is a resolution of singularities).
This implies that the action of $\beta$ on the kernel of the cycle map
$$L _mH _l(X)\rightarrow H _l(X)$$
is zero for any $m\geq n-r-2$. This means that
$$N(Ker(L^mH^*(X)\rightarrow H^*(X)))=0$$
for $0\leq m\leq r+2$

Point a) was proved in \cite{P}.
\end{proof}

\begin{corollary}

Let $X$ be a smooth projective variety such that rational equivalence coincides with homological equivalence in $CH _*(X)\tensor \mathbb{Q}$ in degrees less or equal than $r$. Then the algebraic equivalence coincides with homological equivalence in $CH _*(X)\tensor \mathbb{Q}$ in degrees less or equal than $r+1$ (\cite{P})  and in degrees greater or equal than $n-r-2$.
\end{corollary}
\begin{remark}
\label{r10}
 It is conjectured \cite{Par} that for $X$ a smooth complete intersection in $\mathbb{P}^n$ of multi-degree $d _1\geq d _2\geq ..\geq d_s$ we have 
$$ CH _l(X) _\mathbb{Q}\simeq \mathbb{Q}$$ 
for any $l\leq k-1$ where $k=[\frac{n-\sum _{i=2,s}d _i}{d _1}]$, the integer part of the rational number.

In particular this would imply that $X$ has small Chow group of dimension $k-1$. Supposing this conjecture and using Theorem \ref{t3} we conclude that in our case we have
$$Griff^r(X) _\mathbb{Q}=0$$ 
for any $r\geq n-k$ and $r\leq k+1$ and moreover 
$$K^{r,*} _\mathbb{Q}=0$$
for the same range of indexes.

J. Lewis proved the statements for the Griffiths groups of a generic hypersurface in \cite{Le} (without using the above mentioned conjecture).
\end{remark}
It is known \cite{Par} that the conjecture from Remark \ref{r10} is valid for generic cubic fivefold and sixfold. The next proposition studies the Lawson homology groups of such cubics.

\begin{proposition}
Let $X$ be a smooth generic cubic of dimension $d=5$ or 6. Then Suslin's conjecture is valid for $X$.
\end{proposition}
\begin{proof}
In \cite{Par} it is proved that a generic smooth cubic of dimension $d\geq 5$ has $$CH _0(X)\simeq CH _1(X)\simeq \mathbb{Z}.$$ This implies that there is a decomposition 
$$N\Delta =\alpha _0+\alpha _1+\beta$$
with $\alpha _i$ supported on $V _i\times W _{d-i}$ and $\beta$ supported on $X\times \Gamma _{d-2}$, cycles of codimension $d$ in $X\times X$. As in Proposition \ref{t3} we get the equality $N\Delta _*=\beta _*$ on the kernel and the cokernel of the cycle map $L _mH _*(X)\rightarrow H _*(X)$ for any $m$.

Suppose $d=5$. Because the action of $\beta$ is zero for $m\geq 3$ and the action of $\beta^t$ is zero for $m\leq 5-3=2$ we obtain that the Lawson homology of a generic smooth cubic fivefold is isomorphic with singular homology up to torsion, i.e
$$L _*H _*(X)\tensor \mathbb{Q}\simeq H _*(X)\tensor \mathbb{Q}.$$
Using Corollary \ref{c1} we obtain Suslin's conjecture for generic smooth cubic fivefold.

Suppose now $d=6$. Then the action of $\beta$ is zero for $m\geq 4=dim(\Gamma _{d-2})$ and the action of $\beta^t$ is zero for $m\leq 6-4=2=codim(\Gamma _{d-2})$. We remark that the action of $\alpha _i$ is still zero for any $m\geq 1$. As above, we conclude that for any generic smooth cubic sixfold 
    $$ L _mH _*(X)\tensor \mathbb{Q}\simeq H _*(X)\tensor \mathbb{Q}$$ 
for any $m\geq 4$ and any $m\leq 2$. The action of $\beta$ on $L _3H _l(X)$ factors through $L _3H _l(\Gamma _4)\rightarrow H _l(\Gamma _4)$ which is an isomorphism for any $l\geq 7$ and a monomorphism for $l=6$. It implies that the cycle map $L _3H _6(X)\tensor \mathbb{Q} \rightarrow H _6(X)\tensor Q$ is injective and that $L _3H _l(X)\tensor \mathbb{Q}\simeq H _l(X)\tensor \mathbb{Q}$ for any $l\geq 7$. Using now Corollary \ref{c1} we conclude that Suslin's conjecture holds for any generic smooth cubic sixfold.
\end{proof}
\begin{theorem} 
\label{u}
Let $X$ be a smooth projective variety. If $CH^* _\mathbb{Q}(X)\simeq H^* _\mathbb{Q}(X)$ then $L _*H _*(X) _\mathbb{Q}\simeq H _*(X) _\mathbb{Q}$.

 In particular Suslin's conjecture is valid for such $X$.

\end{theorem} 

\begin{proof}
   Our condition on $X$ gives us the following decomposition of the diagonal 
$$N\Delta =\alpha _0+\alpha _1+...+\alpha _d$$
where each $\alpha _i$ is supported on $V _i\times W _{d-i}$. We remark that $dim(V _i)+dim(W _{d-i})=dimX=d$. Using Corollary \ref{l2} we can conclude that 
$$Nx=N\Delta _*(x)=0$$
for any $x$ in the kernel and in the cokernel of the cycle map
$$L _*H _*(X)\rightarrow H _*(X).$$
This implies that 
$$L _mH _l(X)\tensor {\mathbb{Q}}\simeq H _l(X)\tensor {\mathbb{Q}}.$$
We conclude now Suslin's conjecture for $X$ by using Corollary \ref{c1}.
\end{proof}
We notice that the same techniques used in the Theorem \ref{t3} and in Theorem \ref{u} give us the following proposition (which was already known (see \cite{Vo})).
\begin{proposition} 
Let $X$ be a smooth projective complex variety such that the cycle class maps
$$
cl: CH _l(X) _Q\rightarrow H _{2l}(X,Q)
$$
are injective for $l\leq k$. Then $H^{p,q}(X)=0$ for 

a)$p\neq q$, $p+q$ even and $q\leq k$.

b)$\mid p-q\mid >1$, $p+q$ odd and $q\leq k$.
\end{proposition}

\begin{proof}
The vanishing of the above Hodge numbers come from the equalities in the coniveau filtration generated by the decomposition of the diagonal and from the fact that coniveau filtration with complex coefficients is included in the Hodge filtration.
\end{proof}
\section{The case of projective linear varieties}
 Totaro \cite{To} and Jannsen \cite{Ja} gave the definition of a linear variety (see also \cite{RJ}).
\begin{definition}
A complex variety is called 0-linear if it is either empty set or isomorphic to any affine space $\mathbb{A}^n _\mathbb{C}$. Let $n>0$. A complex variety $Z$ is $n$-linear if there is a triple $(U,X,Y)$ of complex varieties so that $Y\subset X$ is a closed immersion with $U$ its complement; $Y$ and one of the varieties $U$ or $X$ is $(n-1)$-linear and $Z$ is the other member in ${U,X}$. We say that $Z$ is linear if it is $n$-linear for some $n\geq 0$.
\end{definition}
 Among examples of linear varieties are toric varieties, flag varieties (\cite{To},\cite{Ja},\cite{RJ}). R. Joshua \cite{RJ} and B. Totaro\cite{To} proved the following $K\ddot{u}nneth$ formula for projective linear varieties:
\begin{theorem}
\label{tl}
(\cite{To}, \cite{RJ}) Let $X$ a projective smooth linear variety of dimension $d$. Then 
$$CH^*(X)\tensor CH^*(X)\simeq CH^*(X\times X).$$
In particular there is a decomposition of the diagonal cycle $\Delta \in CH^d(X\times X)$ of the form
$$\Delta=\sum \alpha _i\times \beta _i$$
with $\alpha _i$, $\beta _i\in CH^*(X)$ algebraic cycles with $dim(\alpha _i)+dim(\beta _i)=d$.
\end{theorem}
Using Corollary \ref{l2} and Theorem \ref{tl} we can conclude that the action of $\Delta$ is zero on the kernel and the cokernel of the cycle map $L _mH _*(X)\rightarrow H _*(X)$. This implies the following proposition:
\begin{proposition}
Let $X$ a projective smooth linear variety. Then $$L _*H _*(X)\simeq H _*(X).$$
In particular we have $$K^{sst} _i(X)\simeq ku^{-i}(X^{an})$$
for any $i\geq 0$.
\end{proposition}
 The above proposition was first proved in \cite{FHW} by other methods.

%The next proposition is a direct result of the methods that we employed above i%n the text. 
%\begin{proposition}
%Let $X$ be a smooth projective variety of dimension $n\geq 3$ and suppose $CH _%k(X) _Q\hookrightarrow H _{2k}(X) _Q$ for any $k\leq n-3$. Then $L _qH _r\simeq %H _r$ for any $r\geq 2n-2$ and any $q$. In particular $X$ fulfills F-M conject%ur
%e.\end{proposition} 
\bibliographystyle{amsplain}
\bibliography{mircea}
%\begin{thebibliography}{99}
%\bibitem{M}
%C.Peters Lawson homology for varieties with small Chow groups...
%\bibitem{S}
%C.Voisin Hodge theory and complex algebraic geometry II 
%\bibitem{GH}
%J.D. Lewis Three lectures on the Hodge conjecture
%\bibitem{BS}
%Bloch-Srinavas Remarks on correspondences and algebraic cycles
%\bibitem{L}

%E.Friedlander Algebraic cycles, Chow varieties, and Lawson homology
%\bibitem{W}
%M.Walker Morphic Abel-Jacobi map
%\bibitem{B}
%Bloch Lectures on algebraic cycles
%\bibitem{SV}
%Voisin-Soule Torsion cohomology classes and algebraic cycles on complex manifol%ds
%\bibitem{V}
%Voisin On integral Hodge classes on uniruled or Calabi-Yau threefolds
%\bibitem{SV}
%Suslin-Voevodski Bloch-Kato conjecture and motivic cohomology with finite coeff%icients
%\bibitem{IP}
%Iskovskikh-Prokhorov Fano varieties
%\bibitem{FM}
%Friedlander-Mazur Correspondence Homomorphisms for Singular Varieties
%\bibitem{AC}
%Albano-Colino On the Griffiths group of the cubic sevenfold
%\bibitem{PW}
%Pedrini-Weibel The higher K-theory of complex varieties
%\bibitem{FL}
%Friedlander-Lawson A theory of algebraic cocycles
%\bibitem{F}
%Friedlander Bloch-Ogus Properties for Topological Cycle Theory
%\end{thebibliography}
\end{document}